\documentclass[11pt]{article}

\usepackage{graphicx,times,xcolor}

\usepackage[breaklinks=true, colorlinks=true, pdfstartview=FitV, linkcolor=blue, citecolor=blue, urlcolor=red]{hyperref}

\usepackage{amsmath}
\usepackage{amsfonts}
\usepackage{amsthm}

\usepackage{mathabx}

\allowdisplaybreaks

\renewcommand{\Vec}[1]{\mbox{\boldmath$#1$}}

\begin{document}

\title{On mean elements in artificial satellite theory\thanks{Submitted to Celestial Mechanics and Dynamical Astronomy}}
\author{Martin Lara}
\date{\color{red}February 13, 2023}
\maketitle

\begin{abstract}
The merits of a perturbation theory based on a mean to osculating transformation that is pure periodic in the fast angle are investigated. The exact separation of the purely short-period effects of the perturbed Keplerian dynamics from the long-period mean frequencies is achieved by a non-canonical transformation, which, therefore, cannot be computed by Hamiltonian methods. For this case, the evolution of the mean elements strictly adheres to the average behavior of the osculating orbit. However, due to the inescapable truncation of perturbation solutions, the fact that this theory confines the long-period variations of the semimajor axis into the mean variation equations, how tiny they may be, can have adverse effects in the accuracy of long-term semi-analytic propagations based on it.
%While tiny, due to the inescapable truncation of a perturbation solution, these effects may dominate the errors of a long-term semi-analytical propagation based on a pure periodic mean to osculating transformation.
\end{abstract}

%\newpage

\section{Introduction}

Mean elements are useful in a variety of aerospace engineering tasks, like maneuver design or fast orbit propagation. They are roughly understood as the result of an averaging process
%\footnote{\color{red}Borrar \dots We quote from \cite{McClain1977} ``The mean elements are defined operationally as the solution to the averaged equations of motion. Consequently, the exact definition of a particular set of mean elements depends on the interval over which the equations of motion are averaged.''} 
whose aim is to isolate the long trend variations of the osculating orbit, which are modulated with short period fluctuations. While the averaging ensues from a particular transformation of variables it happens that different transformations can be chosen. Therefore, the ``mean elements'' terminology is loosely used with the meaning of elements that are free from short-period effects, a vague definition that may cause confusion in the implementation and use of available perturbation solutions in the literature ---refer to the caveats in \cite{Vallado2013}, pp.~616, 690, and 695.
\par

Ideally, the transformation from mean to osculating elements should be pure periodic in the mean anomaly, so that the mean elements comprise all the secular and long-period variations of the osculating orbit. But this is not often the case, and the mean to osculating transformation of different satellite theories may involve non-periodic and long-period terms in addition to the short-period terms. This is, in particular, the case of common analytical solutions in closed form of the eccentricity, in which some long-period terms remain hidden in the mean to osculating transformation due to the implicit dependence of the true anomaly on the mean anomaly and the eccentricity through Kepler's equation \cite{Kozai1962}. Because of that, the mean orbit will depart from the average osculating orbit in long period displacements of small amplitude.
\par

The lack of agreement between mean and average dynamics is not a concern in usual analytical or semi-analytical orbit theories, a case in which the osculating solution is obtained after recovering the periodic effects from the mean to osculating transformation. This is the case of usual solutions computed by canonical perturbation theory \cite{Brouwer1959,Kozai1959,CoffeyNealSegermanTravisano1995,Lara2021IAC}. However, the disagreement may be unwanted in those cases in which the information provided by the mean elements alone is crucial for the problem at hand, as may happen in some space geodesy applications \cite{WagnerDouglasWilliamson1974,MetrisExertier1995} or when high accuracy is needed in preliminary orbit and maneuver design \cite{BhatFrauenholzCannell1990,Guinn1991}. Conversely, the reported aim of some non-canonical perturbation methods is to achieve the removal of only the short-period fluctuations from the mean orbit (see \cite{McClain1977,DanielsonNetaEarly1994}, and references therein). In this way, it is expected that the mean to osculating transformation is free from long-period terms, and, therefore, the resulting mean elements accurately represent the average behavior of the osculating orbit \cite{MetrisExertier1995}. 
\par

Taking these facts into account is of primordial importance when comparing canonical and non-canonical solutions \cite{LaraSanJuanFolcikCefola2011,LaraSanJuanLopezOchoaCefola2014}. Moreover, special attention must be paid to the correct computation of the short-period elimination when Hamiltonian simplification techniques are involved in the computation of the canonical solution \cite{LaraSanJuanLopezOchoa2013b,LaraSanJuanLopezOchoa2013c}. On top on that some confusion is added since canonical perturbation methods have been reported to be particularizations of non-canonical methods for the appropriate choice of the arbitrary integration functions of the slow variables arising in both different methods \cite{Morrison1966}. But this connection between canonical and non-canonical perturbation theories simply recalls the obvious: that canonical transformations are very particular instances of the wider field of transformation theory. Remarkably, the mean to osculating transformation that produces the exact separation between short- and long-period variations has been recognized as non-canonical for perturbed Keplerian motion \cite{MetrisExertier1995}. More precisely, the exact separation of terms that are pure periodic in the mean anomaly only can be guaranteed up to the first order of the small parameter when using Hamiltonian methods \cite{FerrazMello1999}.
\par

Regarding artificial satellite theory, this issue is mainly relevant in the computation of second order effects of the dominant zonal harmonic of the second degree, hereafter noted $J_2$. Indeed, for Earth-like bodies second order effects of $J_2$ are comparable to those produced by the higher order harmonics of the gravitational potential, and, therefore, are routinely incorporated into operational software \cite{CoffeyAlfriend1984,CoffeyNealSegermanTravisano1995,FolcikCefola2012,LaraSanJuanHautesserres2018}. 
\par

The generalized method of averaging (see \cite{Nayfeh2004} p.~168 and ff., for instance) is a suitable, non-canonical option when approaching semi-analytically the solution of non-trivial perturbation models, which may include both Hamiltonian and non-Hamiltonian perturbations \cite{McClain1977}. Rather, we focus on the simple $J_2$ perturbation and resort to the versatile method of Lie transforms to compute explicit expressions of the mean variations %up to the third order of $J_2$, 
as well as their concomitant periodic corrections. Beyond its original aim of dealing with perturbed Hamiltonian flows \cite{Deprit1969}, the Lie transforms method generally applies to systems of ordinary differential equations with minor modifications of the original formulation, and is readily implemented by means of efficient recursive algorithms \cite{Kamel1970,Henrard1970Lie}. Most notably, the method allows for the explicit computation of the \emph{inverse}, osculating to mean transformation, in addition to the \emph{direct}, mean to osculating one. The former, we will see, is fundamental for understanding the consequences of choosing a mean to osculating transformation that is pure periodic in the mean anomaly in order to obtain the mean frequencies of the motion.
\par

To clearly illustrate the non-canonical nature of a mean to osculating transformation that is pure periodic in the mean anomaly, we strip the formulation of non-essential complexities stemming from the closed form integration of implicit functions of the mean anomaly. That is, we base our exposition on a $J_2$-type potential in which the implicit dependence of the $J_2$-problem on the mean anomaly is avoided by resorting to the usual expansions of the elliptic motion \cite{BrouwerClemence1961}. Furthermore, we truncate these expansions to the lower orders of the eccentricity to shorten the length of printed expressions. %{\color{red}Still, for reference, we also provide the mean variations of orbital elements in closed form for the actual $J_2$ problem.} %  which clearly illustrate that non-canonical and canonical perturbation solutions, while equivalent in the computation of osculating elements, differ in the second order of the mean variation of the mean anomaly.
\par

The toy model is assumed exact in the computation of the mean variations so that the mandatory analytical tests carried out to ensure the correctness of the theory may provide immediate insight without need of making additional expansions. In particular, at each order of the perturbation approach we checked that the successive composition of the mean to osculating and the osculating to mean analytic transformations yields formally the identity up to the truncation order of the theory. We also checked that replacing the transformation equations of the short-period elimination in the original, osculating variation equations, leads to mean variation equations that match term for term the mean frequencies obtained with the perturbation approach up to the truncation order of the perturbation solution. These tests and model also serve us to illustrate the paradox that the inverse of a non-canonical transformation that is pure periodic in the fast angle may bear long-period as well as non-periodic terms.
\par

Moreover, the toy model permits the computation of high orders of the perturbation approach with relatively little computational effort. We need to compute them to show that, for a mean to osculating transformation that is pure periodic in the mean anomaly, the mean variation of the semimajor axis ceases to vanish at the third order. This is in clear contrast with the canonical transformation approach, in which making the mean anomaly ignorable in the Hamiltonian definitely results in constant mean semimajor axis at any order. Because of that, we may anticipate the inadequacy of a perturbation theory based on a pure periodic mean to osculating transformation for long-term propagation. Indeed, because the more important source of errors of a perturbation solution is related to the truncation of the mean frequencies, the neglected long-period effects in the variation of the semimajor axis unavoidably introduce long-period errors that will become dominant in the long term. Tests carried out to check the accuracy of perturbation solutions based on different choices of the mean to osculating transformation confirm this conjecture.
\par

For the same reasons of simplicity and insight, we constructed the perturbation solution in Keplerian orbital elements, which are singular for circular orbits as well as for equatorial orbits. The implementation of alternative non-singular solutions useful for operational purposes can be approached analogously. %{\color{red}Obtaining the solution in a different set of variables, as for instance non-singular ones, does not need re-implementing the perturbation approach, and is just a matter of reformulating the perturbation solution ---as it happens with the Hamiltonian case too \cite{Lyddane1963,DepritRom1970}. However, while this kind of reformulation is certainly possible, we will see that, in general, the generator of a pure periodic mean to osculating transformation in some set of variables does not necessarily yield a transformation of the same character in a different set of variables. In consequence, the non-canonical perturbation theory should be reconstructed in the new variables from scratch if we want to preserve some desired feature ---like the pure periodic character of the mean to osculating transformation.}
%{\color{red}Whether or not this last case corresponds to a canonical transformation remains to be investigated}.
\par

The paper is organized as follows. For completeness, the basic equations of the Lie transforms method are summarized in Section \ref{s:Lie}. Then, the toy model is formulated in usual orbital elements in Section \ref{s:toy}. Also in this section we discuss the two basic instances in which the computation of the $m$-th order terms of the mean variations can be decoupled from the computation of the mean to osculating transformation of the same order, and illustrate the paradox that the inverse of a pure periodic transformation may yield long period effects of higher order than the first. The perturbation solution based on a mean to osculating transformation that is pure periodic in the mean anomaly is described in Section \ref{s:pureperi}, showing that the appearance of long-period terms in the variation of the mean semimajor axis starts at the third order. Section \ref{s:nullC} deals with the case of a perturbation solution based on the pure periodic character of the vectorial generator of the Lie transforms method, which, like in the Hamiltonian case, results in null mean variation of the semimajor axis at any order. Finally, selected propagations that illustrate the main features of each perturbation approach are presented in Section \ref{s:tests}.

\section{The Lie transforms method for vectorial flows} \label{s:Lie}

In this section we present the basic formulation that should allow interested readers to reproduce our computations. Full details in the derivation of the algorithms can be consulted in the original reference \cite{Kamel1970} as well as in the reformulation in \cite{Henrard1970Lie} to which we rather adhere.
\par

A Lie transformation $x_j=x_j(y_k,\epsilon)$, with $j$ and $k$ integers ranging form $1$ to $l$, is defined as the solution of the differential system
\begin{equation} \label{lt}
\frac{\mathrm{d}{x}_j}{\mathrm{d}\epsilon}=W_j(x_k,\epsilon), \qquad j,k=1,2,\dots,l,
\end{equation}
where $\epsilon$ is the independent variable and the $W_j$ denote the $l$ components of some vectorial generating function, for the initial conditions
\begin{equation} \label{iicc}
y_j=x_j(y_k,0).
\end{equation}

%\subsection{Lie transform of a function}
When such kind of transformation is applied to an analytical function $F=F(x_j,\epsilon)$ given by its Taylor series expansion
\begin{equation} \label{fold}
F=\sum_{m\ge0}\frac{\epsilon^m}{m!}\,\left.\frac{\partial^mF}{\partial\epsilon^m}\right|_{\epsilon=0}
=\sum_{m\ge0}\frac{\epsilon^m}{m!}\,F_{m,0}(x_j),
\end{equation}
the transformed function $F^*=F(x_j(y_k,\epsilon),\epsilon)$ can be directly obtained as a Taylor series in the new variables
\begin{equation} \label{fnew} %F^*\equiv
F^*=\sum_{m\ge0}\frac{\epsilon^m}{m!}\,\left.\frac{\mathrm{d}^mF}{\mathrm{d}\epsilon^m}\right|_{\epsilon=0}
=\sum_{m\ge0}\frac{\epsilon^m}{m!}\,F_{0,m}(y_k),
\end{equation}
by the recursive computation of the coefficients $F_{0,m}$ from Deprit's triangle
\begin{equation} \label{deprittriangle}
F_{m,q+1}=F_{m+1,q}+\sum_{i=0}^m{m\choose{i}}\mathcal{L}_{i+1}(F_{m-i,q}), 
\end{equation}
in which $\mathcal{L}$ denotes the linear, scalar operator
\begin{equation} \label{scalarH}
% \mathcal{L}_{i+1}(F_{m-i,q})=\sum_{k=1}^{l}\frac{\partial{F}_{m-i,q}}{\partial{x}_k}{W}_{k,i+1},
\mathcal{L}_{m}(\psi)=\sum_{k=1}^{l}\frac{\partial\psi}{\partial{x}_k}{W}_{k,m},
\end{equation}
and ${W}_{k,m}$ are the coefficients of the Taylor series expansion of the vectorial generator. Namely,
\begin{equation} \label{ltW}
W_j=\sum_{m\ge0}\frac{\epsilon^m}{m!}\,W_{j,m+1}(x_k), \qquad j,k=1,2,\dots,l.
\end{equation}
%{\color{red}
%Note that $\epsilon$ acts as the independent variable of the Lie transformation and also as the small parameter of the Taylor series expansions \cite{Hori1966,Deprit1969}. This duality roughly anticipates that the validity of such kinds of expansions may need to be carefully discussed in longer time intervals than $t=1/\epsilon$ \cite{BreakwellVagners1970}.}
\par

%\subsection{Perturbation solution of a differential system by Lie transforms}
On the other hand, the formulation of a differential system, say
\begin{equation} \label{vectorflow}
\frac{\mathrm{d}x_j}{\mathrm{d}t}=X_j(x_k,\epsilon)\equiv\sum_{m\ge0}\frac{\epsilon^m}{m!}\Phi_{j,m,0}(x_k), \quad
\Phi_{j,m,0}=\left.\frac{\partial^m{X}_j}{\partial\epsilon^m}\right|_{\epsilon=0},
\end{equation}
in the new variables $y_j=y_j(x_k,\epsilon)$ needs to take the Jacobian of the transformation into account. Namely,
\begin{equation} %\label{vectorflownew}
\frac{\mathrm{d}y_j}{\mathrm{d}t}
=\sum_{k=1}^l\frac{\partial{y}_j}{\partial{x}_k}X_j(x_k,\epsilon)=Y_j(x_k,\epsilon),
\end{equation}
whose series expansion in the new variables is
\begin{equation} \label{vectorflownew}
\frac{\mathrm{d}y_j}{\mathrm{d}t}=Y_j(x_k(y_i,\epsilon),\epsilon)
\equiv\sum_{m\ge0}\frac{\epsilon^m}{m!}\Phi_{j,0,m}(y_k), \quad
\Phi_{j,0,m}=\left.\frac{\mathrm{d}^m{Y}_j}{\mathrm{d}\epsilon^m}\right|_{\epsilon=0}.
\end{equation}
Now, the coefficients $\Phi_{j,0,n}$ are obtained directly in the new variables from the new recursion
\begin{equation} \label{vrecurrence}
\Phi_{j,m,q+1}=\Phi_{j,m+1,q}+\sum_{i\ge0}^{m} {m\choose{i}}\mathcal{L}^*_{j,i+1}(\Phi_{k,m-i,q}), \quad
j,k=1,2,\dots,l,
%\sum_{k=1}^l\left(\frac{\partial\Phi_{j,m-i,q}}{\partial{x}_k}W_{k,i+1}-\frac{\partial{W}_{j,i+1}}{\partial{x}_k}\Phi_{k,m-i,q}\right),
\end{equation}
%where $\Vec{\Phi}_{m-i,q}=(\Phi_{1,m-i,q},\Phi_{2,m-i,q},\dots\Phi_{l,m-i,q})$, 
which still adheres to the structure of Deprit's triangle but with a different, vectorial operator under the summations. To wit,
%\[
%\mathcal{L}^*_{j,i+1}(\Phi_{k,m-i,q})=\sum_{k=1}^l\left(\frac{\partial\Phi_{j,m-i,q}}{\partial{x}_k}W_{k,i+1}-\frac{\partial{W}_{j,i+1}}{\partial{x}_k}\Phi_{k,m-i,q}\right)
%\]
\begin{equation} \label{vecop}
\mathcal{L}^*_{j,i}(\Phi_k)=\sum_{k=1}^l\left(\frac{\partial\Phi_j}{\partial{x}_k}W_{k,i}-\frac{\partial{W}_{j,i}}{\partial{x}_k}\Phi_k\right), \quad
j,k=1,2,\dots,l.
\end{equation}
\par

%\subsection{Perturbations by Lie transforms}
In a perturbation approach, we search for a generator that transforms the vectorial flow (\ref{vectorflow}) in some desired, simplified form given by (\ref{vectorflownew}) ---commonly with one of the variables removed up to a given order of the small parameter $\epsilon$. For this task, the first order of recursion (\ref{vrecurrence}) yields
\begin{equation} \label{homo1}
\mathcal{L}^*_{j,1}(\Phi_{k,0,0})=\Phi_{j,0,1}-\tilde\Phi_{j,0,1},
\end{equation}
where $\tilde\Phi_{j,0,1}=\Phi_{j,1,0}$. After choosing the terms $\Phi_{j,0,1}$ in accordance with our simplification criterion, Eq.~(\ref{homo1}) is solved for the $W_{j,1}$ from the system of partial differential equations obtained replacing $\mathcal{L}^*_{j,1}$ from Eq.~(\ref{vecop}). Analogously, at the second order, repeated iterations of Eq.~(\ref{vrecurrence}) yield
\begin{equation} \label{homo2}
\mathcal{L}^*_{j,2}(\Phi_{k,0,0})=\Phi_{j,0,2}-\tilde\Phi_{j,0,2},
\end{equation}
where $\tilde\Phi_{j,0,2}=\Phi_{j,2,0}+\mathcal{L}^*_{j,1}(\Phi_{k,1,0})+\mathcal{L}^*_{j,1}(\Phi_{k,0,1})$ depends only on terms of the first order, and hence is known. Then, after choosing $\Phi_{j,0,2}$ at our convenience, we solve Eq.~(\ref{homo2}) for $W_{j,2}$.
\par

That is, in general, Eq.~(\ref{vrecurrence}) is rearranged as the homological equation
\begin{equation} \label{homological}
\mathcal{L}^*_{j,m}(\Phi_{k,0,0})
%\equiv
%\sum_{k=1}^l\left(\frac{\partial\Phi_{j,0,0}}{\partial{x}_k}W_{k,m}-\frac{\partial{W}_{j,m}}{\partial{x}_k}\Phi_{k,0,0}\right)
=\Phi_{j,0,m}-\tilde\Phi_{j,0,m}, \qquad j=1,\dots,l,
\end{equation}
in which the terms $\tilde\Phi_{j,0,m}$ are known from previous computations, whereas those $\Phi_{j,0,m}$ are chosen in agreement with some desired objective. Therefore, everything becomes known in the homological equation save for the $m$-components of the vectorial generating function (\ref{ltW}), which are then solved from the corresponding partial differential system.
\par

The direct transformation
\begin{equation} \label{dirtra}
x_j=x_j(y_k,\epsilon)\equiv\sum_{m\ge0}\frac{\epsilon^m}{m!}x_{j,0,m}(y_k), \qquad j,k=1,\dots,l,
\end{equation}
is derived from the vectorial generating function $W$ as a particular instance of the scalar recursion (\ref{deprittriangle}) for the functions $x_j=\sum_{m\ge0}(\epsilon^m/m!)\,x_{j,m,0}(x_k)$, in which $x_{j,0,0}=x_j$ and $x_{j,m,0}=0$ for $m\ge1$.
\par

The inverse transformation $y_j=y_j(x_k,\epsilon)$ is analogously obtained from the vectorial generating function $V_j=-W_j$, which in turn must be written in the $y_k$ variables. Because $W$ is itself a vectorial flow, the reformulation of $W=W(x_k(y_k,\epsilon),\epsilon)$ as a power series in the $y_k$ variables is done based on Eq.~(\ref{vrecurrence}), which is readily implemented replacing $B_{j,i,0}=W_{j,i+1}(x_k)$. Once the terms $B_{j,0,i}$ are computed recursively from Eq.~(\ref{vrecurrence}), the inverse transformation
\begin{equation} \label{inverse}
V_j=\sum_{m\ge0}\frac{\epsilon^m}{m!}\,V_{j,m+1}(y_k),
\end{equation}
is obtained by making $V_{j,i+1}=V_{j,i+1}(y_k)\equiv-B_{j,0,i}(x_k)|_{\epsilon=0}$. In particular, $V_1=V_1(y_k)\equiv-W_1(x_k)|_{\epsilon=0}$, from which it follows that the first order terms of direct and inverse transformations are formally opposite ---yet they are evaluated in distinct sets of variables. On the other hand, making $m=q=0$ in Eq.~(\ref{vrecurrence}) immediately yields $B_{j,0,1}=B_{j,1,0}=W_{j,2}$, and hence $V_2=V_2(y_k)\equiv-W_2(x_k)|_{\epsilon=0}$. Remark that this trivial equivalence between direct and inverse generating function terms no longer applies beyond the second order of $\epsilon$.
\par

\section{A $J_2$-type perturbed Keplerian model} \label{s:toy}

The $J_2$ nonspherical term exerts the dominant gravitatonial perturbation of the Keplerian motion on Earth orbiting artificial satellites. This is the reason why this perturbation model is known as the ``main problem'' in the theory of artificial satellites \cite{Brouwer1959}. For this problem, the only disturbance of the two body potential $\mathcal{V}=-\mu/r$ is due to the term
\begin{equation} \label{VJ2}
\mathcal{D}=-\frac{\mu}{r}J_2\frac{R_{\Earth}^2}{r^2}\frac{1}{4}\left[
2-3\sin^2I+3\sin^2I\cos2(f+\omega)\right],
\end{equation}
where $r=a(1-e^2)/(1+e\cos{f})$ is the distance from the center of mass of the attracting body, $a$, $e$, $I$, $\Omega$, $\omega$, and $M$ are usual Keplerian elements, denoting semimajor axis, eccentricity, inclination, right ascension of the ascending node, argument of the periapsis, and mean anomaly, respectively, whereas the true anomaly $f$ is an implicit function of $M$ and involves the solution of the Kepler equation. For a given body, the disturbing potential (\ref{VJ2}) is particularized by the physical parameters $\mu$, $R_{\Earth}$, and $J_2$, standing for the gravitational parameter, the equatorial radius, and the oblateness coefficient, respectively.
\par

The fact that the mean anomaly $M$, which is the fast angle to be removed in the computation of the mean variations, appears implicitly in Eq.~(\ref{VJ2}) as a function of the true anomaly $f$, introduces additional complications in the solution of the integrals involved in the perturbation approach if the closed form is wanted. Therefore, in order to focus on the relevant facts of the solution in mean elements we rely on the usual expansions of the elliptic motion \cite{BrouwerClemence1961}, which, besides, are truncated to the lower orders of the eccentricity. In this way we not only sidestep many of the difficulties arising in the computation of the closed form solution, but also shorten higher order expressions stemming from the perturbation solution to a manageable extent.
\par

%\subsection{The toy model}

Therefore, we replace the disturbing potential of the main problem in Eq.~(\ref{VJ2}), by the substitute
\begin{align} \nonumber
\mathcal{T}= &\; \frac{\mu}{2a}J_2\frac{R_{\Earth}^2}{a^2}\frac{1}{4}\Big[ 2(3s^2-2)(1+3e\cos{M})
+ 3es^2\cos(M+2\omega) \\ \label{VJ2toy}
& -6s^2\cos(2M+2\omega) -21es^2\cos(3M+2\omega) \Big],
\end{align}
that depends explicitly on the mean anomaly, in which we abbreviate $s\equiv\sin{I}$. Then, the variations of the osculating orbital variables are cast in the form of Eq.~(\ref{vectorflow}), where $\epsilon=J_2$, and $\mathrm{d}x_j/\mathrm{d}t$ denote the time variations of $a$, $e$, $I$, $\Omega$, $\omega$, and $M$, for $j=1,\dots6$, respectively. In particular, %$\Phi_{j,0,0}=0$ for $j<6$, $\Phi_{6,0,0}=n$, and
\begin{align} \label{Fj00}
\Phi_{j,0,0}= &\; 0, \qquad j<6,  \\ \label{lp0}
\Phi_{6,0,0}= &\; n  \\ \nonumber
\Phi_{1,1,0}= &\; an\frac{R_{\Earth}^2}{a^2}\frac{3}{4}\big[ 
(6s^2-4)e\sin{M} +es^2\sin(M+2\omega) -4s^2\sin(2M+2\omega) \\ \label{ap}
&  -21es^2\sin(3M+2\omega) \big], \\ \nonumber
\Phi_{2,1,0}= &\; n\frac{R_{\Earth}^2}{a^2}\frac{\eta}{e}\frac{3}{8}\big[ 
e\eta(6s^2-4)\sin{M} +e(\eta-2)s^2\sin(M+2\omega) \\ \label{ep}
& -4(\eta-1)s^2\sin(2M+2\omega) -7e(3\eta-2)s^2\sin(3M+2\omega) \big], \\ \label{ip}
\Phi_{3,1,0}= &\; n\frac{R_{\Earth}^2}{a^2}\frac{cs}{\eta}
\frac{3}{4}\big[ e\sin(M+2\omega) -2\sin(2M+2\omega) -7e\sin(3M+2\omega)\big], \\ \nonumber
\Phi_{4,1,0}= &-n\frac{R_{\Earth}^2}{a^2}\frac{c}{\eta}
\frac{3}{4}\big[ 2 +6e\cos{M} +e\cos(M+2\omega) -2\cos(2M+2\omega) \\ \label{hp}
& -7e\cos(3M+2\omega)\big], \\ \nonumber
\Phi_{5,1,0}= &-n\frac{R_{\Earth}^2}{a^2}\frac{1}{e\eta}\frac{3}{8}\Big\{ 
4e(s^2-1)\big[1-\cos(2M+2\omega)\big] \\ \nonumber
& +\big[e^2(6s^2-8)+6s^2-4\big]\cos{M} +\big[e^2(s^2-2)+s^2\big]\cos(M+2\omega) \\ \label{gp}
& -7\big[e^2(s^2-2)+s^2\big]\cos(3M+2\omega) \Big\}, \\ \nonumber
\Phi_{6,1,0}= &-n\frac{R_{\Earth}^2}{a^2}\frac{1}{e}\frac{3}{8}\Big\{ 4e\big[(3s^2-2-3s^2\cos(2M+2\omega)\big]
+(7e^2-1) \\ \label{lp}
&\times\big[(6s^2-4)\cos{M}+ s^2\cos(M+2\omega) -7s^2\cos(3M+2\omega) \big] \Big\}.
\end{align}
where $n=(\mu/a^3)^{1/2}$ denotes the Keplerian mean motion, and we also abbreviate $c\equiv\cos{I}$. Terms $\Phi_{j,m,0}$ vanish for $m\ge2$. The appearance of the eccentricity in denominators of the first order terms of the variations of the argument of the periapsis and the mean anomaly given in Eqs.~(\ref{gp}) and (\ref{lp}), respectively, may cause trouble in the integration of low eccentricity orbits. However, because the toy model is chosen just for illustrative purposes this is not a concern, and we simply will take care in our numerical experiments of choosing initial conditions far away enough from problematic cases. Certainly, a perturbation solution intended for operational purposes should rather be implemented in some set of non-singular variables.
\par

\subsection{Perturbed Keplerian motion. The homological equation}

Remarkably, for perturbed Keplerian motion the homological equation of the Lie transforms method can be solved by indefinite integration. More precisely, for the $m$-th term, Eq.~(\ref{homological}) turns into
\begin{align} \label{homomJ2}
W_{j,m}= & \;
\frac{1}{n}\int\big(\tilde\Phi_{j,0,m}-\Phi_{j,0,m}\big)\mathrm{d}M, \qquad j=1,\dots5, \\ \label{homo6J2}
W_{6,m}= & \;\frac{1}{n}\int\Big(\tilde\Phi_{6,0,m}-\Phi_{6,0,m}+\frac{\partial n}{\partial{a}}W_{1,m}\Big)\mathrm{d}M,
\end{align}
where $\partial{n}/\partial{a}=-3n/(2a)$, as follows form the definition of the mean motion of a Keplerian flow. An attentive reader will have noticed the similarities between Eqs.~(\ref{homomJ2})--(\ref{homo6J2}) and Eqs.~(55)--(58) in \cite{Hori1971}. %In general, Eq.~(\ref{homo6J2}) cannot be solved until Eq.~(\ref{homomJ2}) is integrated for $j=1$, but we will see in the presented aplications that this is not always the case.
\par

It is important to note that $\Phi_{j,0,m}$ should be chosen in such a way that it cancels the average of the known, tilde terms over the mean anomaly, namely, $\langle\tilde\Phi_{j,0,m}\rangle_M$, in this way preventing the appearance of secular terms in the solution of the homological equation for $j=1,\dots5$. On the other hand, the choice of such $\Phi_{6,0,m}$ that cancels the terms $\langle\tilde\Phi_{6,0,m}-\frac{3}{2}(n/a)W_{1,m}\rangle_M$, thus avoiding the undesired secular terms, relies on the previous computation of $W_{1,m}$. More precisely, we only need to know $\langle{W}_{1,m}\rangle_M=C_{1,m}(a,e,I,\Omega,\omega,-)$, which can be chosen in advance in two relevant cases due to the freedom provided by the arbitrary integration function that arises in the integration of Eq.~(\ref{homomJ2}). In particular, imposing $W_{1,m}$ to be pure periodic in the mean anomaly makes $C_1=0$ trivially, and hence we can choose $\Phi_{6,0,m}=\langle\tilde\Phi_{6,0,m}\rangle_M$ like in the other cases. Alternatively, imposing the mean to osculating transformation of the semimajor axis to be pure periodic in the mean anomaly, allows for the direct computation of $C_{1,m}$ from previous orders of the perturbation solution.
\par

For the last case, the mean to osculating transformation in Eq.~(\ref{dirtra}) of the semimajor axis $a=a'+J_2a_{0,1}+\frac{1}{2}J_2^2a_{0,2}+\dots$ is recursively computed from Eq.~(\ref{deprittriangle}) as follows.  At the first order $a_{0,1}=\mathcal{L}_{1}(a)=W_{1,1}$, as follows from the definition of the scalar operator in Eq.~(\ref{scalarH}). Then, 
\begin{equation} \label{C11}
C_{1,1}=0,
\end{equation}
and hence $\Phi_{6,0,1}=\langle\Phi_{6,1,0}\rangle_M$. At the second order, $a_{0,2}=a_{1,1}+\mathcal{L}_{1}(a_{0,1})$, where $a_{1,1}=W_{1,2}$, as follows from Eq.~(\ref{deprittriangle}). Hence $W_{1,2}=a_{0,2}-\mathcal{L}_{1}(a_{0,1})$, where, for a pure periodic mean to osculating transformation $\langle{a}_{0,2}\rangle_M=0$. Therefore, the computation of $C_{1,2}=\langle{W}_{1,2}\rangle_M$ yields
\begin{equation} \label{C12}
C_{1,2}=-\langle\mathcal{L}_{1}(a_{0,1})\rangle_M
=-\sum_{k=1}^6\Big\langle\frac{\partial{a}_{0,1}}{\partial{x}_k}W_{k,1}\Big\rangle_M,
\end{equation}
that only relies on first order terms. Then, $\Phi_{6,0,2}=\langle\tilde\Phi_{6,0,2}\rangle_M-\frac{3}{2}(n/a)C_{1,2}$. Analogously, at the third order, Eq.~(\ref{deprittriangle}) leads to the sequence
\[
\begin{array}{rl}
a_{2,1}= & W_{1,3}   \\
a_{1,2}= & \mathcal{L}_{2}(a_{0,1}) +\mathcal{L}_{1}(a_{1,1}) +a_{2,1}  \\
a_{0,3}= & \mathcal{L}_{1}(a_{0,2}) +a_{1,2}
\end{array}
\]
from which $W_{1,3}=a_{0,3}-\mathcal{L}_{1}(a_{0,2})-\mathcal{L}_{2}(a_{0,1})-\mathcal{L}_{1}(a_{1,1})$. Therefore, the requirement that $\langle{a}_{0,3}\rangle_M=0$ yields
\begin{equation} \label{C13}
C_{1,3}=-\langle\mathcal{L}_{1}(a_{0,2})\rangle_M -\langle\mathcal{L}_{2}(a_{0,1})\rangle_M -\langle\mathcal{L}_{1}(a_{1,1})\rangle_M,
\end{equation}
and the long-period terms $\langle\tilde\Phi_{6,0,3}\rangle_M-\frac{3}{2}(n/a)C_{1,3}$ that should be cancelled by $\Phi_{6,0,3}$ are computed with the only knowledge of second order terms. And so on.
\par

In a paradox, the choice of arbitrary functions that make the mean to osculating transformation pure periodic may prevent the inverse transformation from having the same nature. Using vectors we denote $\Vec{x}=\Vec{x}'+\epsilon\Vec{x}_{0,1}(\Vec{x}')+\frac{1}{2}\epsilon^2\Vec{x}_{0,2}(\Vec{x}')+\mathcal{O}(\epsilon^3)$ the mean-to-osculating transformation, and $
\Vec{x}'=\Vec{x}+\epsilon\Vec{x}'_{0,1}(\Vec{x})+\frac{1}{2}\epsilon^2\Vec{x}_{0,2}'(\Vec{x})+\mathcal{O}(\epsilon^3)$ its inverse. Neglecting terms $\mathcal{O}(\epsilon^3)$ and higher, their sequential composition yields
\begin{equation}
\Vec{x}=\Big[\Vec{x}+\epsilon\Vec{x}'_{0,1}(\Vec{x})+\frac{1}{2}\epsilon^2\Vec{x}_{0,2}'(\Vec{x})\Big]
+\epsilon\Vec{x}_{0,1}(\Vec{x}+\epsilon\Vec{x}'_{0,1}(\Vec{x}))
+\frac{1}{2}\epsilon^2\Vec{x}_{0,2}(\Vec{x}),
\end{equation}
where the term $\epsilon\Vec{x}_{0,1}$ still needs to be expanded. That is,
\begin{equation}
\Vec{0}= \epsilon\left[\Vec{x}'_{0,1}(\Vec{x})+\Vec{x}_{0,1}(\Vec{x})\right]+\frac{1}{2}\epsilon^2\left[
\Vec{x}_{0,2}'(\Vec{x}) 
+2\Vec{x}'_{0,1}(\Vec{x})\frac{\partial\Vec{x}_{0,1}(\Vec{x})}{\partial\Vec{x}}
+\Vec{x}_{0,2}(\Vec{x})\right],
\end{equation}
in which $\Vec{x}'_{0,1}(\Vec{x})$ and $\Vec{x}_{0,1}(\Vec{x})$ cancel each other out because they are just opposite. Hence,
\begin{equation} \label{x02p}
\Vec{x}_{0,2}'=-2\Vec{x}_{0,1}(\Vec{x})\frac{\partial\Vec{x}_{0,1}(\Vec{x})}{\partial\Vec{x}}-\Vec{x}_{0,2}(\Vec{x}),
\end{equation}
where $\Vec{x}_{0,2}$ is pure periodic in $M$ by our choice of the arbitrary integration functions. However, the product $\Vec{x}_{0,1}\partial\Vec{x}_{0,1}/\partial\Vec{x}$ may give rise ---and it certainly does for the $J_2$ problem--- to constant and long-period terms in spite of each factor is pure periodic in the mean anomaly.%\footnote{A simple illustration is provided by $\cos^2\alpha$}
\par

Details on the construction of perturbation solutions based in these two noteworthy cases are provided in the two following Sections.

\section{Theory 1} \label{s:pureperi}

In this Section we obtain the perturbation solution based on a mean to osculating transformation that is pure periodic in the mean anomaly, which we label as ``Theory 1''.
\par

At the first order $n=0$ and $q=0$ in Eq.~(\ref{vrecurrence}), and hence $\tilde\Phi_{j,0,1}=\Phi_{j,1,0}$. Besides, $C_{1,1}=0$, as follows from Eq.~(\ref{C11}), and we choose $\Phi_{j,0,1}=\langle\tilde\Phi_{j,0,1}\rangle_M$. We obtain $\Phi_{1,0,1}=\Phi_{2,0,1}=\Phi_{3,0,1}=0$, and
\begin{equation}
\Phi_{4,0,1}=-n\frac{R_\Earth^2}{a^2}\frac{3c}{2\eta}, \quad
\Phi_{5,0,1}=-c\Phi_{4,0,1}, \quad
\Phi_{6,0,1}=n\frac{R_\Earth^2}{a^2}\frac{3}{2}(3c^2-1).
\end{equation}
Then, the trivial integration of Eqs.~(\ref{homomJ2})--(\ref{homo6J2}) for $m=1$ provides the first order terms $W_{j,1}$ of the vectorial generating function, each of which will depend on an arbitrary integration function $C_{j,1}\equiv{C}_{j,1}(a,e,I,\Omega,\omega,-)$ save for $C_{1,1}$, whose value has already been fixed.
\par

The mean to osculating transformation is of the form of Eq.~(\ref{dirtra}), and is computed using recursion (\ref{deprittriangle}) in which $F_{m,0}=0$ for $m\ge1$ and $F$ is successively replaced by each orbital element $x_j\in(a,e,I,\Omega,\omega,M)$. In these equations, the arbitrariness of the functions $C_{j,1}$ is suppressed by imposing that the transformation from mean to osculating elements be pure periodic in the mean anomaly, which yields $C_{j,1}=0$. Then, we obtain
\begin{align} \nonumber
W_{1,1}=a_{0,1}=& -a\frac{R_{\Earth}^2}{a^2}\frac{3}{4}\big[ e(6s^2-4)\cos{M} +es^2\cos(M+2\omega) -2s^2
\\ \label{a01}
&  \times\cos(2M+2\omega) -7es^2\cos(3M+2\omega)  \big], \\ \nonumber
W_{2,1}=e_{0,1}=&-\frac{R_{\Earth}^2}{a^2}\frac{1}{8}\Big[6\eta^2(3s^2-2)\cos{M} +3\eta(\eta-2)s^2\cos(M+2\omega) \\ \label{e01}
& +\frac{6e\eta s^2}{1+\eta}\cos(2M+2\omega) -7\eta(3\eta-2)s^2\cos(3M+2\omega) \Big], \\ \nonumber
W_{3,1}=I_{0,1}=& -\frac{R_{\Earth}^2}{a^2}\frac{cs}{4\eta}\big[ 3e\cos(M+2\omega) -3\cos(2M+2\omega) \\ \label{i01}
&-7e\cos(3M+2\omega)\big], \\ \nonumber
W_{4,1}=\Omega_{0,1}= &-\frac{R_{\Earth}^2}{a^2}\frac{c}{4\eta}\big[ 18e\sin{M} +3e\sin(M+2\omega) -3\sin(2M+2\omega)
\\ \label{h01}
& -7e\sin(3M+2\omega)\big], \\ \nonumber
W_{5,1}=\omega_{0,1}=& -\frac{R_{\Earth}^2}{a^2}\frac{1}{e\eta}\frac{1}{8}\Big\{  6[e^2(3s^2-4) +3s^2-2]\sin{M} \\ \nonumber
& +3[e^2(s^2-2)+s^2]\sin(M+2\omega) -6 e(s^2-1) \\ \label{g01}
& \times\sin(2M+2\omega) -7[e^2(s^2-2)+s^2]\sin(3M+2\omega) \Big\}
\\ \nonumber
W_{6,1}=M_{0,1}= &-\frac{R_{\Earth}^2}{a^2}\frac{1}{8e}\Big\{ (4e^2-1)\big[6(3s^2-2)\sin{M} 
+3s^2\sin(M+2\omega) \\ \label{l01}
& -7s^2\sin(3M+2\omega)\big]-9es^2\sin(2M+2\omega) \Big\}.
\end{align}
The reader will forgive the abuse of notation with the aim of reducing the number of printed expressions, for the components of the generator $W_{j,m}$ depend on the osculating orbital variables whereas the periodic corrections are functions of the mean variables. This change of the osculating variables by the mean ones in the corrections, which also affects the mean variations, should be carried out at the end of the whole procedure in order to run the perturbation  theory in a semi-analytical way. On the other hand, the first order terms of the inverse, osculating to mean transformation, are formally opposite to those in Eqs.~(\ref{a01})--(\ref{l01}), but they are evaluated in the original, osculating variables.
\par

At the second order, the solution of Eqs.~(\ref{homo6J2})--(\ref{homomJ2}) needs the previous computation of the tilde terms $\tilde\Phi_{j,0,2}$. Repeated application of Eq.~(\ref{vrecurrence}) yields
\begin{equation} 
\tilde\Phi_{j,0,2}=\Phi_{j,2,0}+\sum_{k=1}^6\left[
\frac{\partial(\Phi_{j,0,1}+\Phi_{j,1,0})}{\partial{x}_k}W_{k,1}
-\frac{\partial{W}_{j,1}}{\partial{x}_k}(\Phi_{k,0,1}+\Phi_{k,1,0})\right],
\end{equation}
whose evaluation only involves straightforward operations. Then, the second order terms of the mean variations of $a$, $e$, $I$, $\Omega$, and $\omega$, are chosen by averaging, to get
\begin{align} \label{F102}
\Phi_{1,0,2}=& \; 0, \\ \nonumber
\Phi_{2,0,2}=&-n\frac{R_{\Earth}^4}{a^4}\frac{9}{16}\frac{s^2e}{1+\eta}\left[11\eta^2(3s^2-2)+3\eta(13s^2-10)+12(s^2-1)\right] \\
& \times\sin2\omega, \\ \nonumber
\Phi_{3,0,2}= &\; n\frac{R_{\Earth}^4}{a^4}\frac{9}{16}\frac{1}{\eta^2}\frac{cse^2}{1+\eta}\left[11\eta^2(3s^2-2)+3\eta(13s^2-10)+12(s^2-1)\right] \\
& \times\sin2\omega, \\ \nonumber
\Phi_{4,0,2}= & n\frac{R_{\Earth}^4}{a^4}\frac{3c}{16\eta^2}\big\{
33\eta^3(43s^2-12)+104\eta^2(s^2-1)-81\eta(15s^2-4) \\ \nonumber
& -116(2s^2-1)
+6\left[11\eta^2(3s^2-1)+3\eta(13s^2-5)+6(2s^2-1)\right] \\
& \times(\eta-1)\cos2\omega 
\big\} ,\\ \nonumber
\Phi_{5,0,2}= & \; n\frac{R_{\Earth}^4}{a^4}\frac{3}{32\eta^2}\big\{
33\eta^3(43s^4-86s^2+16) +52\eta^2(3s^4-6s^2+4) -162\eta \\ \nonumber
& \times(s^2-1)(15s^2-4) -116(s^2-1)(5s^2-2)
+6\big[11\eta^3(3s^4-6s^2  \\ \nonumber
& +2) +\eta^2(9s^4-16s^2+8) -18\eta(s^2-1)(3s^2-1) -6(s^2-1) \\ \label{F502}
&  \times(5s^2-2) \big]\cos2\omega
\big\}
\end{align}
Finally, we compute $C_{1,2}=\langle{W}_{1,2}\rangle_M$ from Eq.~(\ref{C12}), to obtain
\begin{align} \nonumber
C_{1,2}= &\; -a\frac{R_{\Earth}^4}{a^4}\frac{3}{32\eta}
\big\{ 21\eta^3(43s^4-24s^2+8) +104\eta^2s^2(s^2-2) \\ \nonumber
&  -15\eta(45s^4-24s^2+8) -232s^2(s^2-1)  \\ \label{C12pupe}
& +6(\eta-1)\left[7\eta^2(3s^2-2) +\eta(27s^2-22)-12c^2\right]s^2\cos2\omega
\big\},
\end{align}
from which $\Phi_{6,0,2}=\langle\tilde\Phi_{6,0,2}\rangle_M-\frac{3}{2}(n/a)C_{1,2}$, and hence
\begin{align} \nonumber
\Phi_{6,0,2}= & \; n\frac{R_{\Earth}^4}{a^4}\frac{3}{64\eta}\big\{
459\eta^3(43s^4-24s^2+8) +1456\eta^2s^2(s^2-2) -315\eta \\ \nonumber
& \times(45s^4-24s^2+8) -2784s^2(s^2-1)
+6s^2\big[153\eta^3(3s^2-2) \\ \label{F602pupe}
&  +28\eta^2(3s^2-4) -105\eta(3s^2-2) -144(s^2-1)\big]\cos2\omega \big\}.
\end{align}
\par

It follows the solution of Eqs.~(\ref{homomJ2}) and (\ref{homo6J2}) for $m=2$ to obtain the components of the vectorial generator $W_{j,2}$, which again introduce arbitrary integration functions that are particularized for obtaining a pure periodic mean to osculating transformation. In addition to the function $C_{1,2}$ in Eq.~(\ref{C12pupe}), for $j\ge2$ we readily obtain $C_{j,2}=0$.
\par

For instance, for the semimajor axis we obtain
\begin{equation} \label{a02pupe}
a_{0,2}=a\frac{R_{\Earth}^4}{a^4}\frac{3}{32}\frac{1}{1+\eta}\sum_{i=0}^2\sum_{j=1+\delta_{2,i}}^{2i+2}\sum_{k=\delta_{2,i}}^{2+2\delta_{2i,j}}P_{i,j,k}\eta^{k-1}e^{|j-2i|}s^{2i}\cos(2i\omega+jM),
\end{equation}
where $\delta_{i,j}$ denotes the Kronecker delta and the inclination polynomials $P_{i,j,k}$ are presented in Table \ref{t:a02}. That $a_{0,2}$ is pure periodic in the mean anomaly results from the subindex $j\ge1$. On the contrary, the inverse, osculating to mean correction to the semimajor axis takes the form
\begin{equation} \label{a02p}
a_{0,2}'=a_{0,2,\mathrm{long}}'+a\frac{R_{\Earth}^4}{a^4}\frac{3}{16}\sum_{i=0}^2\sum_{j=1}^{2i+2}\sum_{k=1-\delta_{2,i}}^{1+2\delta_{2i,j}}P'_{i,j,k}\eta^{k-1}e^{|j-2i|}s^{2i}\cos(2i\omega+jM),
\end{equation}
with the new inclination polynomials $P'_{i,j,k}$ of Table \ref{t:ap02}, and the long-period terms
\begin{equation} \label{a02plp}
a_{0,2,\mathrm{long}}'=
a\frac{R_{\Earth}^4}{a^4}\frac{3}{16}\sum_{i=0}^1\sum_{k=0}^{3}P'_{i,0,k}\eta^{k-1}s^{2i}\cos2i\omega,
\end{equation}
with coefficients $P'_{i,0,k}$ in the same Table, which clearly show the non pure periodic nature of the osculating to mean transformation announced in Eq.~(\ref{x02p}). Indeed, it is immediate to see that $\langle{a}_{0,2}'\rangle_M=a_{0,2,\mathrm{long}}'$. The same features are indeed observed in the direct and inverse corrections to the other orbital variables.
\par

\par

\begin{table}[htb] % \tabcolsep 4pt
\center
\caption{Inclination polynomials $P_{i,j,k}$ in Eq.~(\protect\ref{a02pupe}).}
\label{t:a02}
\begin{tabular}{@{}lll@{}}
\hline\noalign{\smallskip}
${}_{0,1,0}:-208 s^2 (s^2-1)$ & ${}_{1,2,2}:8 (255 s^2-142)$ & ${}_{2,3,1}:3$ \\
${}_{0,1,1}:-4 (289 s^4-196 s^2+48)$ & ${}_{1,2,3}:-8 (405 s^2-298)$ & ${}_{2,3,2}:-3$ \\
${}_{0,1,2}:-8 (125 s^4-72 s^2+24)$ & ${}_{1,2,4}:-968 (3 s^2-2)$ & ${}_{2,4,1}:-54$ \\
${}_{0,2,0}:224 s^2 (s^2-1)$ & ${}_{1,3,0}:40 (s^2-1)$ & ${}_{2,4,2}:-54$ \\
${}_{0,2,1}:2 (247 s^4-40 s^2-24)$ & ${}_{1,3,1}:2 (281 s^2-194)$ & ${}_{2,4,3}:70$ \\
${}_{0,2,2}:2 (191 s^4+72 s^2-24)$ & ${}_{1,3,2}:162 (3 s^2-2)$ & ${}_{2,4,4}:70$ \\
${}_{1,1,0}:168 (s^2-1)$ & ${}_{1,4,0}:-336 (s^2-1)$ & ${}_{2,5,1}:-63$ \\
${}_{1,1,1}:6 (109 s^2-82)$ & ${}_{1,4,1}:-84 (7 s^2-6)$ & ${}_{2,5,2}:-77$ \\
${}_{1,1,2}:198 (3 s^2-2)$ & ${}_{1,4,2}:-140 (3 s^2-2)$ & ${}_{2,6,1}:-147$  \\
${}_{1,2,0}:720 (s^2-1)$ & ${}_{2,2,1}:-3$ & ${}_{2,6,2}:-147$ \\
${}_{1,2,1}:72 (43 s^2-32)$ & ${}_{2,2,2}:-3$ &  \\
\noalign{\smallskip}\hline
\end{tabular}
\end{table}
\begin{table}[htb] % \tabcolsep 4pt
\caption{Non-null coefficients $P'_{i,j,k}$ in Eqs.~(\protect\ref{a02p}) and (\protect\ref{a02plp}).}
\label{t:ap02}
\begin{tabular}{@{}lll@{}}
\hline\noalign{\smallskip}
${}_{0,0,0}: -232 s^2 (s^2-1)$ & ${}_{0,1,1}: 96 (3 s^4-3s^2+1)$ & ${}_{1,3,1}: -108 (3 s^2-2)$ \\
${}_{0,0,1}: -15 (45 s^4-24 s^2+8)$ & ${}_{0,2,1}: 6 (11 s^4-24s^2+8)$ & ${}_{1,4,1}: -84 (3 s^2-2)$ \\
${}_{0,0,2}: 104 s^2 (s^2-2)$ & ${}_{1,1,0}: 24(s^2-1)$ & ${}_{2,2,1}: 3$ \\
${}_{0,0,3}: 21 (43 s^4-24 s^2+8)$ & ${}_{1,1,1}: -12 (3s^2-2)$ & ${}_{2,3,1}: -12$ \\
${}_{1,0,0}: -72 (s^2-1)$ & ${}_{1,2,0}: -24(s^2-1)$ & ${}_{2,4,1}: -30$ \\
${}_{1,0,1}: -30 (3 s^2-2)$ & ${}_{1,2,1}: -96 (3s^2-2)$ & ${}_{2,4,3}: 42$ \\
${}_{1,0,2}: 12 (3 s^2-4)$ & ${}_{1,2,3}: 72 (3s^2-2)$ & ${}_{2,5,1}: 84$ \\
${}_{1,0,3}: 42 (3 s^2-2)$ & ${}_{1,3,0}: -56(s^2-1)$ & ${}_{2,6,1}: 147$ \\
\noalign{\smallskip}\hline
\end{tabular}
\end{table}
%

%Due to our choice of the arbitrary functions $C_{j,2}$, for $j>1$, t

This apparent inconsistency in the osculating to mean transformation does not corrupt the solution, which, as desired is made of the mean elements resulting from the numerical integration of the mean frequencies, and the analytic pure periodic corrections. Still, this kind of solution may be not the best one when the theory is intended for accurate long-term orbit propagation. To show that we compute the third-order terms of the mean variations. Recalling that $\Phi_{j,k,0}=0$ for $k\ge2$, from repeated iterations of recursion (\ref{vrecurrence}) we obtain
\[
\begin{array}{l}
\Phi_{j,0,3}= \mathcal{L}^*_{j,1}\left(\Phi_{j,0,2}\right)+\Phi_{j,1,2} \\
\Phi_{j,1,2}= \mathcal{L}^*_{j,2}\left(\Phi_{j,0,1}\right)+\mathcal{L}^*_{j,1}\left(\Phi_{j,1,1}\right)+\Phi_{j,2,1} \\
\Phi_{j,2,1}= % L_{W_1}\left(\Phi_{j,2,0}\right)+
\mathcal{L}^*_{j,3}\left(\Phi_{j,0,0}\right)+2\mathcal{L}^*_{j,2}\left(\Phi_{j,1,0}\right)+\Phi_{j,3,0},
\end{array}
\]
where, from the same recursion, $\Phi_{j,1,1}=\Phi_{j,0,2}-\mathcal{L}^*_{j,1}\left(\Phi_{j,0,1}\right)$. Then, the computable terms are $\tilde\Phi_{j,0,3}=\mathcal{L}^*_{j,1}\left(\Phi_{j,0,2}+\Phi_{j,1,1}\right)
+\mathcal{L}^*_{j,2}\left(\Phi_{j,0,1}+2\Phi_{j,1,0}\right)$, and, like before, we compute the third order components of the mean variations as $\Phi_{j,0,3}=\langle\tilde\Phi_{j,0,3}\rangle_M$ for $j<6$, whereas the computation of $\Phi_{6,0,3}=\langle\tilde\Phi_{6,0,3}-\frac{3}{2}(n/a)W_{1,3}\rangle_M$ requires the preliminary computation of $C_{1,3}=\langle{W}_{1,3}\rangle_M$ from Eq.~(\ref{C13}). For the latter we obtain
\begin{align} \nonumber
\Phi_{1,0,3}= &\; an\left(\!\frac{R_{\Earth}}{a}\!\right)^6\frac{c^2s^2e^2}{\eta^2(1+\eta)}\frac{81}{16}\left[
7\eta^2(3s^2-2)+\eta(27 s^2-22)+12(s^2-1)
\right] \\
&\times \sin2\omega.
\end{align}
\par

The non-vanishing character of the mean variation of the semimajor axis is a crucial fact that remains hidden when the perturbation approach is truncated to the lower orders. Indeed, due to the Lyapunov instability of the Keplerian motion, the semimajor axis is the more sensitive element in the propagation of errors. So neglecting these terms of the variation of the mean semimajor axis, how tiny they may be, soon or later will introduce observable errors in the mean elements propagation. This feature of Theory 1 will be clearly illustrated in the semi-analytical propagations presented in Section \ref{s:tests}.
\par

Finally, it is worth mentioning that the formulation of Theory 1 in a different set of variables, as for instance, non-singular ones based on the semi-equinoctial elements $\kappa=e\cos\omega$, $\sigma=e\sin\omega$, while feasible, will loose the pure periodic character of the mean to osculating transformation as soon as we reach the second order. The reason is the same stated above that the product of two pure periodic terms yields, in general, non-periodic terms. Thus, for instance, up to second order terms $\kappa=(e'+\epsilon{e}_{0,1}+\frac{1}{2}\epsilon^2e_{0,2})\cos(\omega'+\epsilon\omega_{0,1}+\frac{1}{2}\epsilon^2\omega_{0,2})
$, which, after expansion and truncation to $\mathcal{O}(\epsilon^3)$ yields
\begin{align} \nonumber
\kappa=& \; e'\cos\omega'
+\epsilon\left(e_{0,1}\cos\omega'-\omega_{0,1}e'\sin\omega'\right)
+\frac{1}{2} \epsilon^2\left(e_{0,2}\cos\omega'-\omega_{0,2}e'\sin\omega'\right) \\
&-\frac{1}{2}\epsilon^2\left(2e_{0,1}\omega_{0,1}\sin\omega'+\omega_{0,1}^2e'\cos\omega'\right),
\end{align}
where the products $e_{0,1}\omega_{0,1}$ and $\omega_{0,1}^2$ will naturally spring long-period terms that will remain after reformulation of the orbital elements into semi-equinoctial variables. Therefore, Theory 1 should be recomputed in the new variables from scratch if we want to preserve the pure periodic character in the fast variable of the mean to osculating transformation.
\par

\section{Theory 2} \label{s:nullC}

Other case in which the successive order of the mean variations can be computed based only on previous orders of the vectorial generator is when we prescribe that the vectorial generating function be pure periodic in the mean anomaly. That is, the arbitrary integration functions $C_{j,m}=\langle{W}_{j,m}\rangle_M$ vanish at any order of the perturbation theory. We label this case as Theory 2.
\par

To the first order, both approaches match and the results in Section \ref{s:pureperi} apply also to this case. At the second order, for $j<6$ the mean frequencies $\Phi_{j,0,2}$ are the same as those given in Eqs.~(\ref{F102})--(\ref{F502}), but the integration function $C_{1,2}$ in Eq.~(\ref{C12pupe}) must be replaced by $C_{1,2}=0$. Hence, 
\begin{align} \nonumber
\Phi_{6,0,2}= & \; n\frac{R_{\Earth}^4}{a^4}\frac{3}{32\eta}\big\{
198\eta^3(43s^4-24s^2+8) +572\eta^2s^2(s^2-2) -135\eta(45s^4 \\ \nonumber
& -24s^2+8)-1044s^2(s^2-1) 
+6s^2\big[ 66\eta^3(3s^2-2) +11\eta^2(3s^2-4) \\ \label{F602C0}
& -45\eta(3s^2-2)-54(s^2-1) \big]\cos2\omega \big\},
\end{align}
whose structure is analogous to Eq.~(\ref{F602pupe}) but with different coefficients. The second order terms of the vectorial generator $W_{j,2}$ are consequently the same as before for $j<6$, but change for $W_{6,2}$. Since they all are involved in the computation of each transformation equation, now all of them are affected of long-period terms, which at this order are exactly the same in the direct and inverse transformations.
\par

In particular, now 
\begin{equation} \label{a02C0}
a_{0,2}=a_{0,2,\mathrm{short}}+a_{0,2,\mathrm{long}},
\end{equation}
where $a_{0,2,\mathrm{short}}$ is the same as Eq.~(\ref{a02pupe}), which only comprises short-period terms, and the long-period terms $a_{0,2,\mathrm{long}}$ are exactly one half of those reported in Eq.~(\ref{a02plp}). Analogously, now
\begin{equation} \label{ap02C0}
a'_{0,2}=a'_{0,2,\mathrm{short}}+a'_{0,2,\mathrm{long}},
\end{equation}
where the short-period terms $a'_{0,2,\mathrm{short}}$ are those obtained from Eq.~(\ref{a02p}) as $a_{0,2}'-a_{0,2,\mathrm{long}}'$, whereas $a'_{0,2,\mathrm{long}}=a_{0,2,\mathrm{long}}$. The same happens to the other variables, in which the long-period terms appearing in the osculating to mean transformations of the previous approach are now distributed in equal parts between the direct and inverse transformation.
\par

Contrary to Theory 1, the current choice of a pure short-periodic vectorial generator makes that now the variation of the mean semimajor axis vanishes at each order. Besides, since there are not special requisites related to the periodic character of the short-period elimination, Theory 2 can be reformulated in a different set of variables without need of approaching the Lie transformations from scratch, as it also happens in the Hamiltonian case \cite{Lyddane1963,DepritRom1970}.
\par

\section{Semi-analytical propagations. A test case } \label{s:tests}

The aim of this section is just to illustrate the different behavior between the two non-canonical perturbation theories described in the previous sections. The actual relevance of the differences observed between both kinds of solutions regarding the implementation of operational software falls out of the scope of the current study. %This kind of investigation would require the use of more accurate perturbation models that the toy model used in the exposition, on the one hand, as well as the formulation of the theory in non-singular variables, on the other hand. 
Therefore, we are satisfied with running our tests for a favorable case in the orbital variables used, which we borrow from \cite{CoffeyAlfriend1984}. In particular, we limit the reported tests to the case of an elliptical orbit with eccentricity $e=0.2$, inclination $I=20^\circ$, and semimajor axis $a=9500$ km. The gravity field parameters used in the propagations are $\mu=398600.4415\;\mathrm{km^3/s^2}$, $R_\Earth=6378.1363$ km, and $J_2=0.001082634$.
\par

First of all, we recall that a semi-analytic propagation is made of three parts. The osculating to mean transformation is applied first to the initial conditions in order to get the initial state in mean variables, as well as, more importantly, to initialize the mean frequencies of the motion. Next, the mean variations are integrated numerically for the mean initial conditions, which is achieved with very long steps. Finally, osculating elements are obtained at each evaluation step from the mean-to-osculating transformation. Because the more important source of errors stems from the mean elements propagation, usual perturbation theories include a higher order of the mean variations than the order of the periodic corrections \cite{Brouwer1959,Kozai1959,DepritRom1970,McClain1977}. On the other hand, the initialization process is crucial to avoiding the abnormal growth of errors in the along-track direction, so the osculating to mean transformation is commonly patched including higher order information related to the semimajor axis conversion \cite{LyddaneCohen1962,BreakwellVagners1970,Lara2020arxiv,Lara2022}. %Early1986
In this fashion, we denote as $m$-th order a theory that comprises the $(m+1)$-th order effects of the mean frequencies and the $m$-th order effects of the transformation equations, but in which the osculating to mean transformation is patched with the $(m+1)$-th order effects of the semimajor axis. The usual calibration of the solution using the energy equation is not used in this case to strictly adhere to the general formulation of a non-Hamiltonian problem.
\par

The first order effects of both distinct perturbation approaches yield identical terms. However, in our notation of a 1st-order theory they disagree in $a'_{0,2}$, cf.~Eqs.~(\ref{a02p})--(\ref{a02plp}) and Eq.~(\ref{ap02C0}), as well as in $\Phi_{6,0,2}$, cf.~Eqs.~(\ref{F602pupe}) and (\ref{F602C0}). These small differences in the formulation produce small but significant changes in the time-history of the errors resulting from each perturbation theory, which are computed with respect to a reference orbit obtained from the numerical integration of the osculating variations provided by Eqs.~(\ref{Fj00})--(\ref{lp}).
\par

\begin{figure}[htbp]
\centerline{
\includegraphics[scale=0.67]{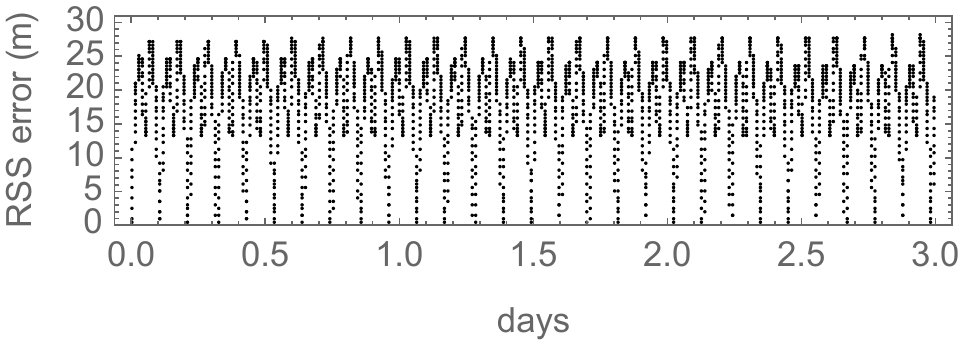}\quad
\includegraphics[scale=0.67]{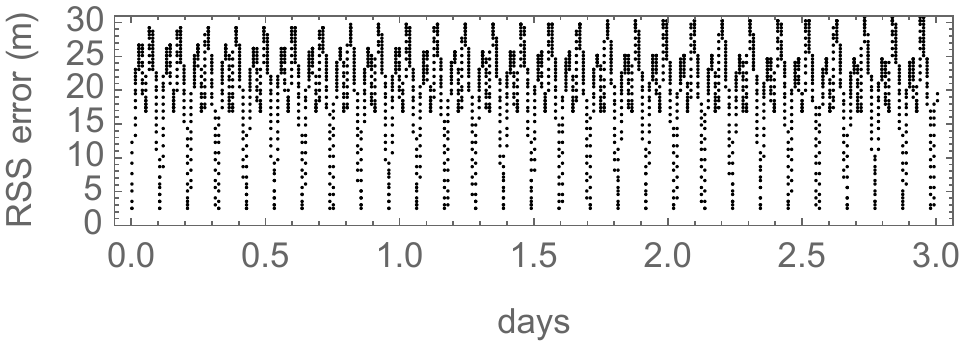}
}
\caption{Test case. Position errors of the 1st-order of Theory 1 (left) and 2 (right).}
\label{f:RSS2nd}
\end{figure}

The time-histories of the root sum square (RSS) of the position errors for a three-day semi-analytic propagation of the test case with each Theory are shown Fig.~\ref{f:RSS2nd}, where we observe that the RSS errors remain of the same order and present the same behavior. This kind of representation is useful in illustrating the accuracy of each theory for orbit prediction purposes, but does not illustrate relevant differences between the two distinct perturbation approaches.
\par

\begin{figure}[htbp]
\centerline{
\includegraphics[scale=0.67]{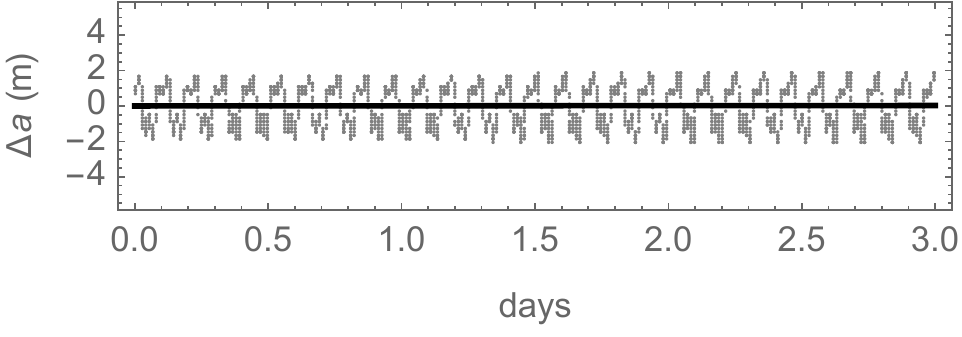} \quad
\includegraphics[scale=0.67]{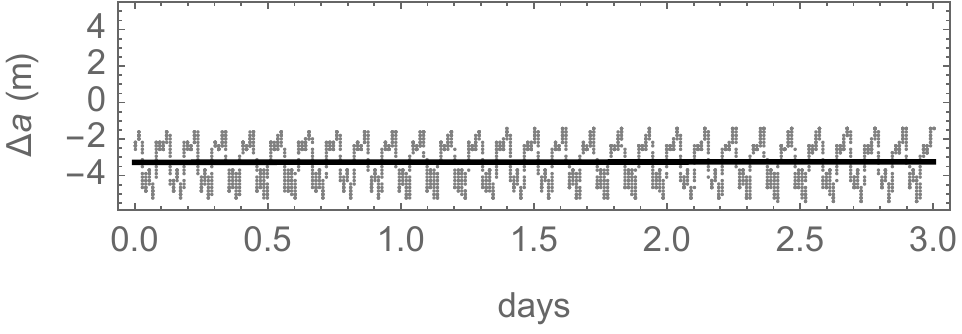}
}
\caption{Test case. Semimajor axis errors of the 1st order of Theory 1 (left) and 2 (right).}
\label{f:a2nd}
\end{figure}

More precisely, Theory 1 yields mean elements that strictly adhere to the average behavior of the osculating orbit, as expected, which is not the case of Theory 2. This is illustrated in Fig.~\ref{f:a2nd} for the semimajor axis, which shows that the propagation errors are of the same order and behave almost identically, with periodic oscillations of the same amplitude. However, the average of the errors, depicted by horizontal black lines in each plot of Fig.~\ref{f:a2nd}, nearly vanishes for Theory 1, averaging to about 1 cm in this example, whereas it is affected of a clear shift from the zero average in the case of Theory 2, of about 3 m in this particular propagation. Differences in the time-histories of the errors of the other orbital elements obtained with each theory, as well as in their respective averages, are of minor nature, and are not presented. In particular, with both theories the errors of the eccentricity average to $\sim10^{-6}$ for the test case, and to values below one tenth of arc second for the orbital angles.
\par

\begin{figure}[htbp]
\centerline{
\includegraphics[scale=0.67]{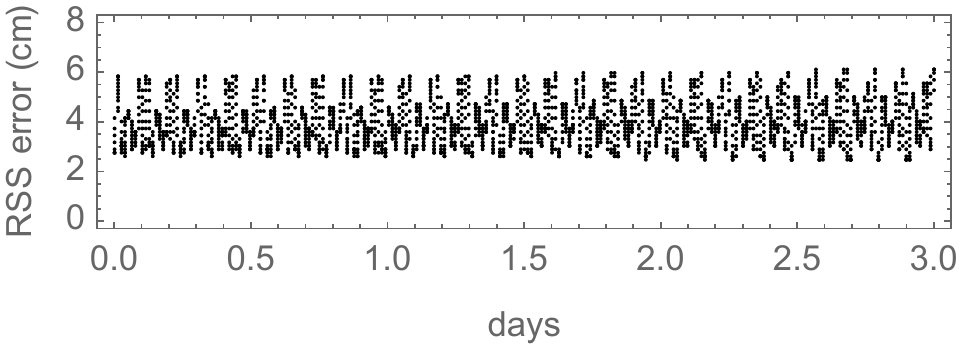} \quad
\includegraphics[scale=0.67]{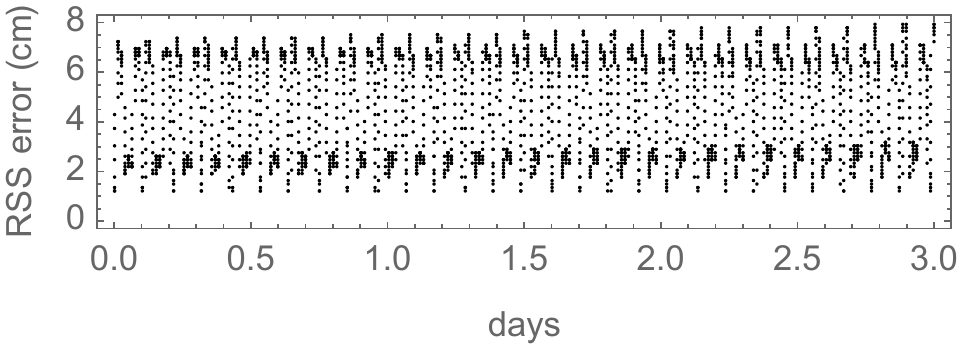} }
\caption{Test case. Position errors of the 2nd-order of Theory 1 (left) and 2 (right).}
\label{f:RSS3rd}
\end{figure}

Taking additional terms clearly improves the accuracy of each perturbation theory, whose position errors reduce from tens of meters to several centimeters when propagating the test case with 2nd-order theories, in agreement with a $\mathcal{O}(J_2)$ improvement. This is shown in Fig.~\ref{f:RSS3rd}, where we also note the lower amplitude of the RSS errors resulting from Theory 1. Nevertheless, in view of the RSS errors remain in both cases with values of the expected order of the truncation of the perturbation theory, we consider that both theories are of comparable accuracy for the test case along the 3-day interval. In particular, in both propagations the errors are $\mathcal{O}(10^{-8})$ relative to distance.
\par

\begin{figure}[htbp]
\centerline{
\includegraphics[scale=0.67]{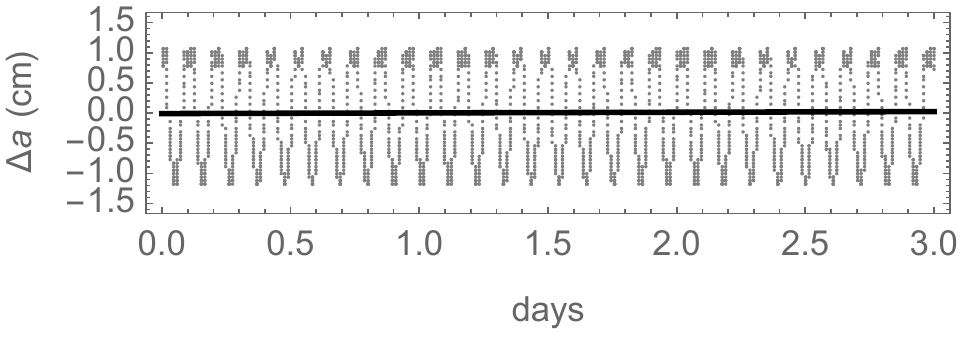} \quad
\includegraphics[scale=0.67]{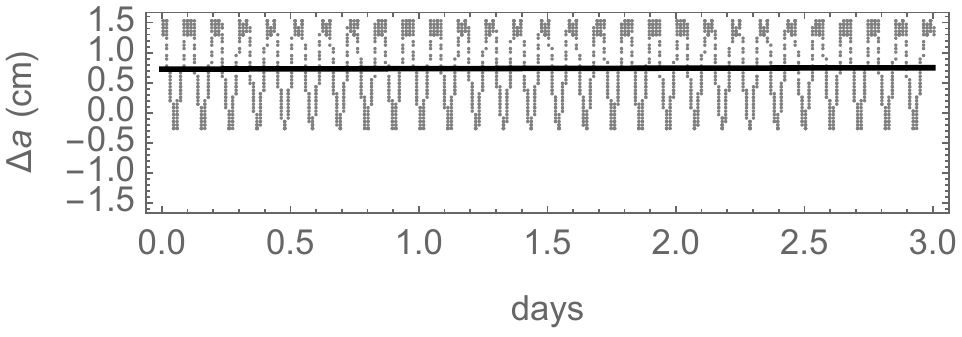} }
\caption{Test case. Semimajor axis errors of the 1st order of Theory 1 (left) and 2 (right).}
\label{f:a3rd}
\end{figure}

Correspondingly, the propagation errors of the semimajor axis reduce their amplitude from the meter to the centimeter level, as depicted in Fig.~\ref{f:a3rd}, also in agreement with the expected improvement provided by an additional order of the perturbation theories. As shown in the right plot of Fig.~\ref{f:a3rd}, the shift from the zero-average in Theory 2 remains, but now reduced to less than 1 cm due to the general improvements of the perturbation solution. The semimajor axis errors stemming from the semi-analytical propagation using Theory 1 fall now below the mm level. Regarding the errors of the other orbital elements, they remain small as well as similar in both cases, like it happened to the lower order solutions.
\par

\begin{figure}[htbp]
\centerline{
\includegraphics[scale=0.67]{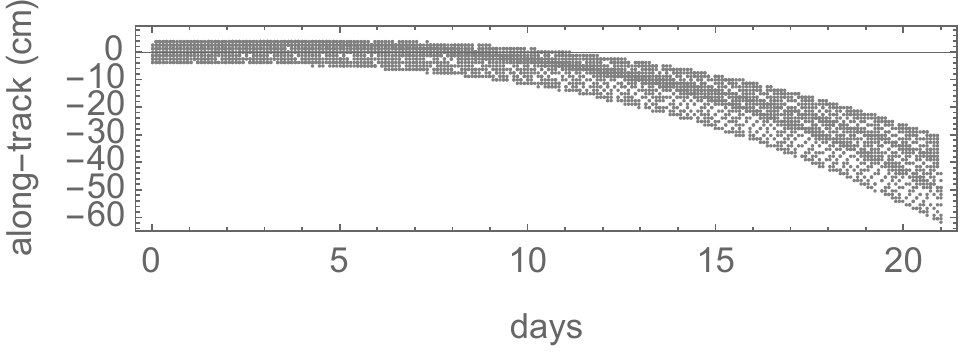} \quad
\includegraphics[scale=0.67]{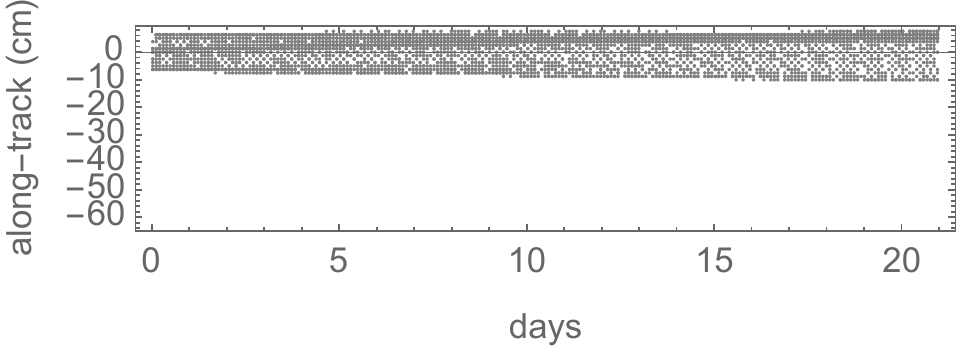} }
\centerline{
\includegraphics[scale=0.67]{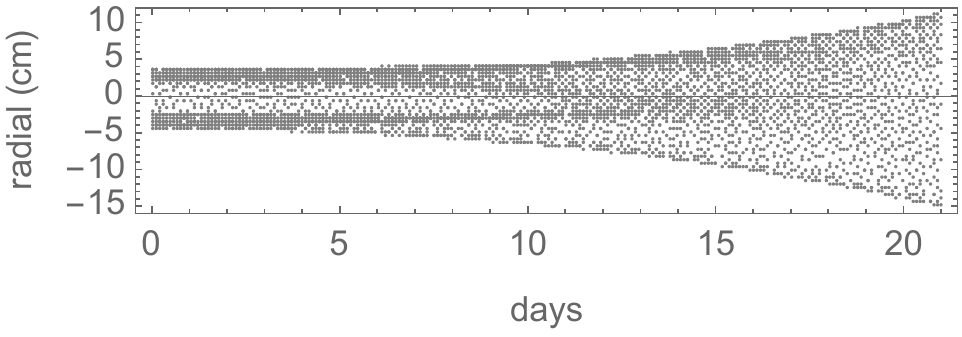} \quad
\includegraphics[scale=0.67]{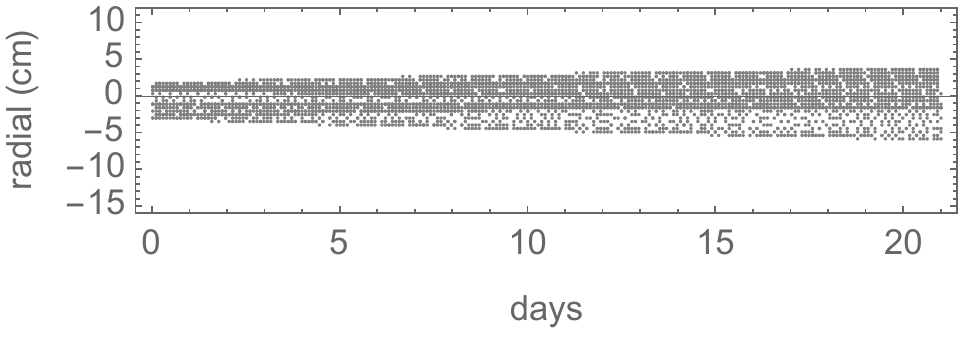} }
\centerline{
\includegraphics[scale=0.67]{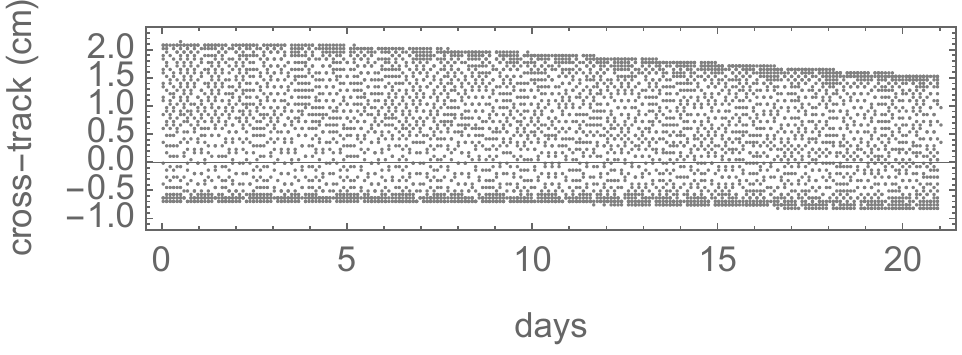} \quad
\includegraphics[scale=0.67]{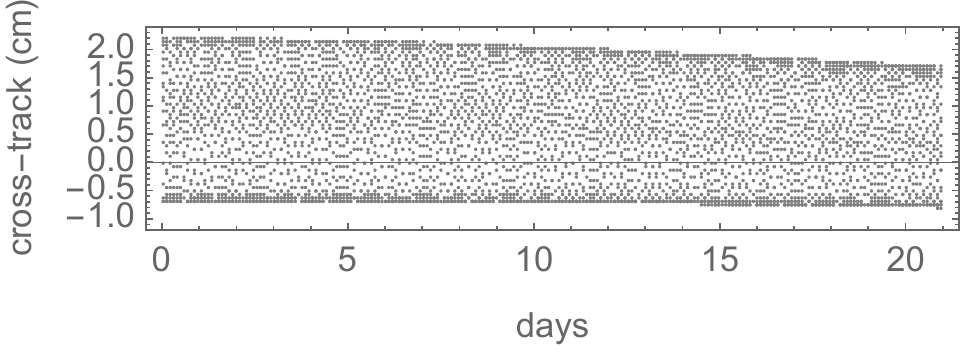} }
\caption{Intrinsic errors of a 3-week semi-analytical propagation of the test case with 2nd-order Theory 1 (left) and 2 (right)}
\label{f:3ia3m2d}
\end{figure}

On the other hand, as already pointed out in the text, the apparent superiority of Theory 1 must be better qualified. For the truncation of the non-vanishing mean variation of the semimajor axis, additional errors mainly affecting the mean motion, and hence the mean anomaly propagation turn this perturbation theory into a less accurate tool for long-term propagation. This is illustrated for the second order versions of the perturbation theories in Fig.~\ref{f:3ia3m2d}, in which the propagation interval of the test case is extended from 3 days to 3 weeks. Now, rather than RSS errors, we display the projections of the position error in the intrinsic (along-track, radial, and cross-track) directions, which better illustrate the nature of the resulting errors. Thus, we clearly see that along-track errors resulting from Theory 1, notably deteriorate passed one week and a half, whereas those of Theory 2 only worsen slightly (top plots of Fig.~\ref{f:3ia3m2d}). This is the expected consequence of the truncation error of the mean motion in the propagation of the mean anomaly, which does not affect Theory 2 as far as the mean variation of the semimajor axis vanishes at any order. Due to the elliptic character of the test orbit, this error in the mean anomaly propagation is also observed in the radial errors, yet with less adverse effects (center plots of Fig.~\ref{f:3ia3m2d}), but barely affects the errors in the cross-track direction (bottom plots of Fig.~\ref{f:3ia3m2d}).
\par

\begin{figure}[htbp]
\centerline{
\includegraphics[scale=0.67]{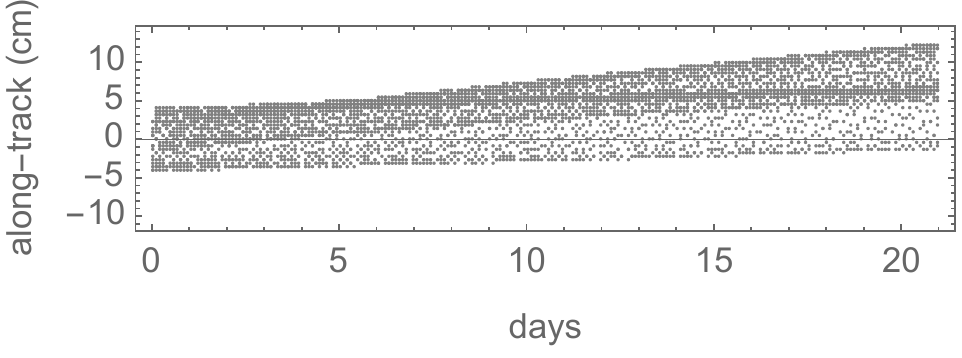} \quad
\includegraphics[scale=0.67]{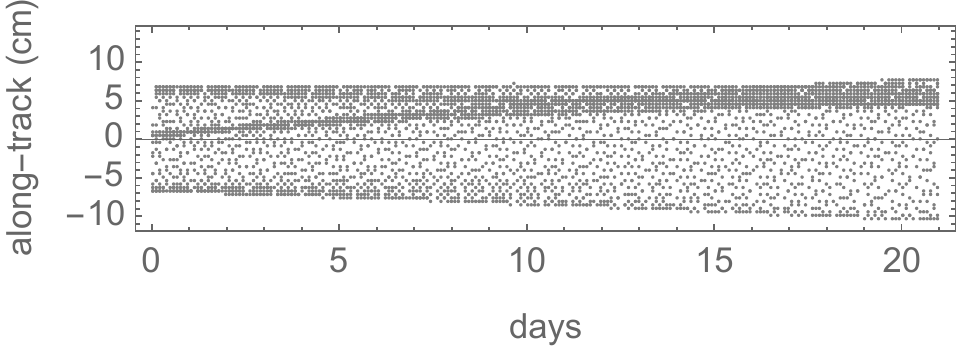} }
\centerline{
\includegraphics[scale=0.67]{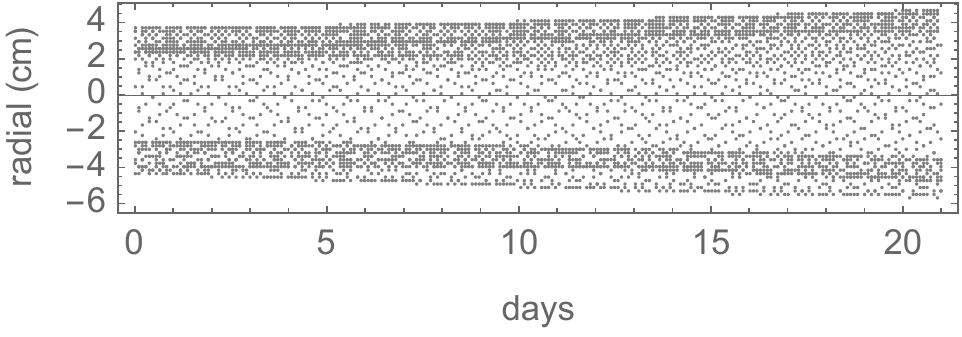} \quad
\includegraphics[scale=0.67]{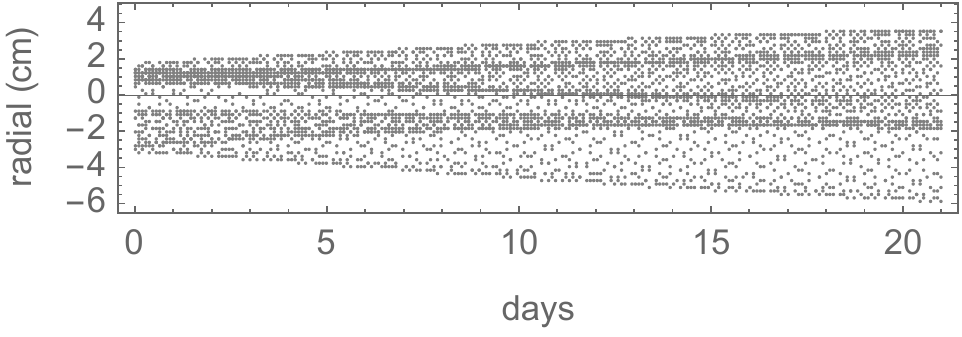} }
\centerline{
\includegraphics[scale=0.67]{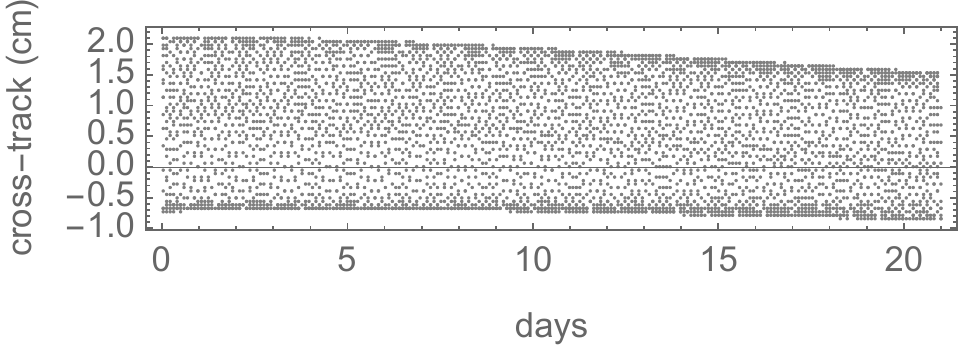} \quad
\includegraphics[scale=0.67]{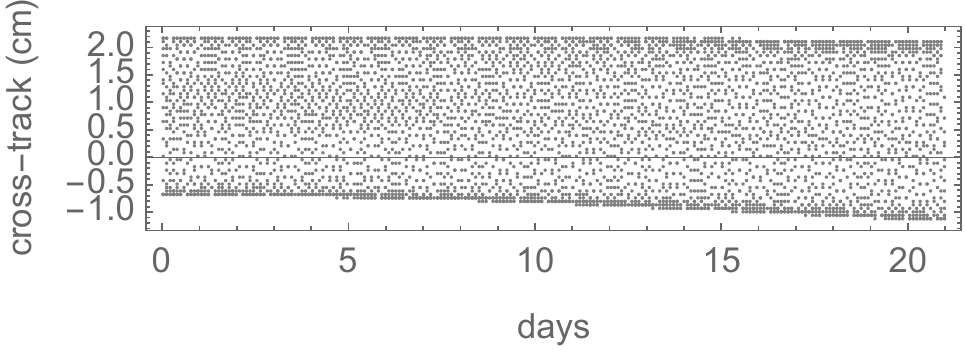} }
\caption{Intrinsic errors of a 3-week semi-analytical propagation of the test case with 2nd-order patched Theory 1 (left) and 2 (right)}
\label{f:3ia4m2d}
\end{figure}

Still, there is no doubt that Theory 1 is certainly correct. To check that, we patched both second order theories with the fourth-order terms of the mean variation of the semimajor axis, which changes nothing in Theory 2 for the vanishing of the mean semimajor axis in this approach, but clearly refines Theory 1. These additional corrections bring both theories again to comparable accuracy for orbit propagation purposes, as shown in Fig.~\ref{f:3ia4m2d}, in which the plots in the right column are the same as corresponding ones in the right column of Fig.~\ref{f:3ia3m2d} but represented in different scales to ease comparison. It is worth to mention that in this last case, the 2nd-order patched Theory 1 provides errors of the semimajor axis that average to just a few hundredths of mm, whereas the average of the semimajor axis error of Theory 2 remains in the cm level in spite of the expanded propagation interval. As before, the errors of the other orbital elements presented by each distinct theory are analogous and very small. On the other hand, the amplitude of the errors remain similar to the previous case because the patched theories still rely only on 2nd-order direct corrections.
\par

Finally, we must mention that, beyond the illustration purposes of this Section, both perturbation solutions of the toy model have been also obtained in non-singular variables ---recomputed for Theory 1, to preserve the pure periodic character of the mean to osculating transformation, and simply reformulated in the case of Theory 2--- in this way allowing us to carry out additional tests for less elliptical orbits that, while non-exhaustive, confirm the reported characteristics of each distinct perturbation approach.

\section*{Conclusions}

For perturbed Keplerian motion, a non-canonical perturbation theory based on a mean to osculating transformation that is pure periodic in the fast angle has the advantage of providing mean elements that strictly adhere to the average evolution of the osculating orbit. On the other hand, extending this kind of perturbation theory to higher orders shows that the mean variation of the semimajor axis only vanishes up to second order of $J_2$ effects. Because the accuracy of a perturbation solution for a given time is essentially related to the truncation order of the mean frequencies, and on account of the Lyapunov instability which is inherent to Keplerian motion, this later fact turns the theory based on the pure periodic mean to osculating transformation less accurate for long-term propagation than classical perturbation approaches due to increasing errors in the along-track direction. Therefore, there is not an always-best perturbation approach, and the choice of the proper kind of perturbation theory must be done depending on the user's particular needs.
\par

For didactic purposes the exposition avoids the use of implicit functions of the mean anomaly and relies on a toy model in common orbital elements, which has been derived from the $J_2$ problem by means of usual expansions of the elliptic motion. Computing an analogous theory in closed form requires to confront the solution of non-trivial integrals. Still, most of the solutions of these integrals are already reported in the literature or can be approached with known techniques of integral calculus. The computation of such kind of perturbation solution is under development, and will be published elsewhere when fully completed and tested. 
\par

%{\color{magenta}
%Perturbations methods for vectorial flows have much wider applicability than canonical perturbation methods, yet at the cost of a non-negligible increase of the computational burden required in the derivation of the perturbation solution. On the other hand, in spite of they can deal effectively with both Hamiltonian and non-Hamiltonian problems, their usefulness in the implementation of simplification techniques analogous to those that are customary in Hamiltonian perturbation methods still remains to be demonstrated. The non-trivial difficulties confronted in this process can be checked by interested readers trying to implement the vectorial elimination of the parallax without using previous knowledge provided by the corresponding Hamiltonian algorithm.
%\par
%}

\subsection*{Acknowledgements}
The research has received funding from the Spanish State Research Agency and the European Regional Development Fund (Projects PID2020- 112576GB-C22 and PID2021-123219OB-I00, AEI/ERDF, EU)

\end{document}